# Numerical Algorithmic Science and Engineering within Computer Science
## Rationale, Foundations, and Organization

John Lawrence Nazareth

March, 2019

**Abstract**

A recalibration is proposed for "numerical analysis" as it arises specifically within the broader, embracing field of modern computer science (CS). This would facilitate research into theoretical and practicable models of real-number computation at the foundations of CS, and it would also advance the instructional objectives of the CS field. Our approach is premised on the key observation that the great "watershed" in numerical computation is much more marked between finite- and infinite-dimensional numerical problems than it is between discrete and continuous numerical problems. A revitalized discipline for numerical computation within modern CS can more accurately be defined as "numerical algorithmic science & engineering (NAS&E), or more compactly, as "numerical algorithmics," its focus being the algorithmic solution of numerical problems that are either *discrete*, or *continuous* over a space of *finite dimension*, or a *combination* of the two. It is the counterpart within modern CS of the numerical analysis discipline, whose primary focus is the algorithmic solution of *continuous*, *infinite-dimensional* numerical problems and their *finite-dimensional approximates*, and whose specialists today have largely been repatriated to departments of mathematics. Our detailed overview of NAS&E from the viewpoints of rationale, foundations, and organization is preceded by a recounting of the role played by numerical analysts in the evolution of academic departments of computer science, in order to provide background for NAS&E and place the newly-emerging discipline within its larger historical context.

## 1. Introduction

In a survey of numerical analysis that appeared in the *Princeton Companion to Mathematics* (Gowers et al. [1]), the renowned numerical analyst, Professor Nick Trefethen of Oxford University, makes the following observation:

"Numerical analysis sprang from mathematics; then it spawned the field of computer science. When universities began to found computer science departments in the 1960s, numerical analysts were often in the lead. Now, two generations later, most of them are to be found in mathematics departments. What happened? A part of the answer is that numerical analysts deal with continuous mathematical problems, whereas computer scientists prefer discrete ones, and it is remarkable how wide a gap that can be."

The repatriation of numerical analysts mentioned above had an important consequence for the mathematics of computation. Motivated by the well-developed, discrete models of computation, which lie at the foundation of computer science and comprise its so-called "grand unified theory" (see, for example, *The Nature of Computation* by Moore and Mertens [2]), the famed mathematician and Fields Medalist winner, Stephen Smale, along with his co-workers, created a counterpart for numerical analysis: real-number, continuous models of computation. In their landmark monograph, *Complexity and Real Computation,* Blum, Cucker, Shub, and Smale [3] presented a model of great generality---abstract machines defined over mathematical rings and fields---and a corresponding theory of computational complexity, in particular, over the real and complex number fields. These "real-number" models serve as a foundation for numerical analysis within mathematics. To foster such activities, Smale et al. also created an umbrella organization known as the Foundations of Computational Mathematics (FoCM) Society [4].

Today, the solution of continuous, infinite-dimensional problems, for example, partial-differential equations, receives the lion's share of attention from numerical analysts. However, such problems are often solved by a reduction to problems of finite dimension, and thus the latter subject also remains an important component of numerical analysis. Numerical computation within computer science (CS), on the other hand, has increasingly focused on solving problems of a *discrete or combinatorial* nature. But it is important to note that continuous problems, *especially problems defined over a real-number space of finite dimension*, also fall centrally within the province of computer science. They arise most commonly in conjunction with discrete numerical problems, good examples being provided by mixed integer programming (linear, nonlinear, and stochastic), network-flow programming, and dynamic programming. The downgrading of numerical analysis as a subfield of computer science noted in the quotation above had a deleterious consequence for CS, because it left largely unfinished the task of building a solid theoretical and practicable foundation for numerical computation within CS *as it relates to real numbers,* one that is *well integrated with the classical, discrete models of computer science*. During the earlier period when numerical analysts had a closer affiliation with computer science, this objective was achieved, but only in a very limited, albeit practically important way, through the development and study of the finite-precision, floating-point model of computation and its associated round-off error analysis.

The primary purpose of this article is to discuss the means whereby this downgrading of real-number computation within CS can be remedied. First, we describe the rationale for explicitly identifying a discipline within computer science, which we term *Numerical Algorithmic Science and Engineering (NAS&E)*, or more compactly, *Numerical Algorithmics*, the counterpart of numerical analysis within mathematics. (Speaking metaphorically, numerical analysis within mathematics can be characterized as seeking to bring "algorithm under the rubric of number" and its NAS&E counterpart within computer science as seeking to bring "number under the rubric of algorithm".) Next, we survey some recently proposed models of real-number computation that can contribute to the task of building a solid theoretical and practical foundation for NAS&E within CS. Finally, we conclude with a brief discussion of the role of NAS&E within the broader CS curriculum, from the standpoints of education and research. Our hope is that this emerging discipline will reoccupy the region within academic departments of computer science that was left vacant following the repatriation of numerical analysts to mathematics (see the above observation of Trefethen quoted from [1]).

We begin, however, with some detailed historical background material that places NAS&E within its broader context. Specifically, we discuss the following: the foundational topic of "number" at the root of mathematics; the corresponding foundational topic of "algorithm and universal computer" at the root of computer science; and two other related historical topics, namely, the key role played by numerical analysts in the setup of university departments of computer science, and the creation of real-number models at the foundation of numerical analysis after its repatriation to mathematics. This material constitutes the first half of our article (Section 2, pgs. 4-11) and it can be perused very quickly or skipped entirely by a reader who prefers to go directly to Sections 3 and 4 (pgs. 12-25), the main thrust of our article as outlined in the preceding paragraph. Afterwards, "GOTO 2", in the language of ancient Fortran, in order to view these proposals from a broader, historical perspective.

## 2. Historical Background

### *2.1 Number*

Few would argue with the observation that mathematics in full-flower as we know it today, both pure and applied, has evolved from the root concept of number. For instance, this is beautifully recounted in *Number:The Language of Science* by Tobias Dantzig [5], a landmark book that was first published in 1930 and then appeared in several subsequent editions. Albert Einstein is said to have endorsed it as follows (italics mine):

"This is beyond doubt the most interesting book on the *evolution of mathematics* that has ever fallen into my hands. If people know how to treasure the truly good, this book will attain a lasting place in the literature of the world. The evolution of mathematical thought from the earliest times to the latest constructions is presented here with admirable consistency and originality and in a wonderfully lively style."

Nowadays every schoolchild learns *number representation* at an early age, along with the basic arithmetic operations on decimal *numerals*. But the *concept of number* itself is far from elementary, a fact highlighted by the British mathematician D.E. Littlewood, a distinguished algebraist, in a chapter titled `Numbers' of his classic, *A Skeleton Key of Mathematics* [6]:

"A necessary preliminary for any proper understanding of mathematics is to have a clear conception of what is meant by number. When dealing with number most people refer to their own past handling of numbers, and this is, usually, not inconsiderable. Familiarity gives confidence in the handling, but not always an insight into the significance. The technique of manipulating numbers is learned by boys and girls at a very tender age when manipulative skill is fairly easily obtained, and when the understanding is very immature. At a later stage, when the faculty of understanding develops, the technique is already fully acquired, so that it is not necessary to give any thought to numbers. To appreciate the significance of numbers it is necessary to go back and reconsider the ground which was covered in childhood. Apart from specialized mathematicians, few people realize that, for example, the number [represented by] 2 can have half a dozen distinct meanings. These differences in meaning are reflected in the logical definitions of number."

Littlewood then proceeds to explain these "differences in meaning" and gives a brief yet masterful exposition of the logical foundations of four basic number systems: cardinal integers, signed integers, rational numbers, and real numbers. In particular, he elaborates on the need for real numbers as follows:

"It is pertinent to enquire why it is necessary to introduce real numbers since these constitute so vast an extension of the rationals, apparently to so little effect, since to every real number, one can obtain a rational approximation to any degree of accuracy. The necessity for the real numbers is illustrated by an important class of theorems called the *existence theorems*. A query often arises, does there exist a number with such and such a property? With rationals the answer is often "no", whereas with real numbers the answer would be "yes". To make sure that a number will always be existent and ready when it is required, the vast extension of rationals to reals is necessary."

The underlying *structure* of the foregoing number systems was subsequently extended, generalized, or relaxed, leading to many other key mathematical concepts, for example, vector spaces, matrix and tensor algebra, groups, rings, and fields, functional analysis, and so on, evolving, over time, into the highly-elaborated mathematics of today.

## 2.2 Algorithm and Universal-Machine

Computer science, in contrast to mathematics, is a much younger discipline, although it too has roots that stretch back to antiquity. Its key foundational concepts are *algorithm* and *universal-machine*, and these foundations were laid during the 1930s---well before the advent of the electronic computer---by a group of mathematical logicians, including Kurt Godel, Alonzo Church, Stephen Kleene, A.A. Markov, Emil Post, and, above all, Alan Turing. (Useful background information on these pioneers and their contributions can be found in Berlinski [7].) Seemingly different formulations of "algorithm" were shown to be to be equivalent, leading to what became known within computer science as the Church-Turing thesis or, more colorfully, as the "*grand unified theory of computation*." For a good overview, see Moore and Mertens [2].

Academic departments of computer science themselves only came into existence within universities during the 1960s and, in rare instances, in the 1950s. (One of the first was created by the great computer pioneer Maurice Wilkes at Cambridge University, where it was known initially as the "mathematical laboratory" and later grew into the university's department of computer science.) In an invaluable collection of articles discussing the underlying philosophy of computer science, one of the founding fathers of this field, Donald Knuth [8], makes the following observation (italics mine):

"My favorite way to describe computer science is to say that it is the study of algorithms...... Perhaps the most significant discovery generated by the advent of computers will turn out to be that algorithms, as *objects of study,* are extraordinarily rich in interesting properties; and, furthermore, that an algorithmic point of view is a useful way to organize information in general. G.E. Forsythe has observed that "the question: `What can be automated?' is one of the most inspiring philosophical and practical questions of contemporary civilization."

Interestingly enough, in the beginning there was some debate as to whether this new discipline should be called *Algorithmics*; see again Knuth [8] who, in turn, attributes the name to Traub [9].

In the public discourse, a "recipe" in a cookbook is often used as an analogue for "algorithm." But, in reality, recipe (say within a soup-cookbook) stands in relation to algorithm in the same way that numeral stands in relation to number. Like number, the concept of "algorithm" is far from elementary: the analogue of an algorithm is, in fact, closer to an *entire* chapter of the soup-cookbook, with different choices of ingredients (inputs to an algorithm) leading, via a sequence of procedural steps, to different soups (outputs). Note, in particular, that an algorithm must always come to a halt and produce its output after a *finite* number of steps. If this is not the case, for example, if on some particular input, its procedural steps enter an infinite loop, then we will use the term *program*. In other words, a program---itself a closer analogue to "recipe"---is not required to produce an answer for each and every given input. (Another frequently used term in this setting is "*computer* program," the concrete realization of an algorithm or program as a finite list of instructions in a computer programming language.) We see that every algorithm is a program within its prescribed model of computation, but every program is *not* necessarily an algorithm. Hence the ubiquitous, so-called halting problem, a key breakthrough of Turing, which can be stated very simply as follows: within a given model of computation, does there exist a particular ("halting") *algorithm* that can examine any *program* whatsoever within the model and determine whether or not that program is an algorithm? The answer, which is premised on the notion of a "universal Turing machine" and utilizes a delightfully simple and elegant (implicit diagonalization) argument, is that such an algorithm *cannot* exist and that the halting problem is therefore *undecidable*, one of the foundational results of the aforementioned "grand unified theory"; see, for example, Stewart [10] or Moore and Mertens [2].

During the 1970s and 80s, the study of the rate-of-convergence and computational complexity of algorithms came to the fore, most notably the breakthroughs of Stephen Cook and Richard Karp on NP-completeness and the identification of the all-encompassing *P=NP* problem of theoretical computer science. For an early, yet comprehensive overview, see Garey and Johnson [11].

The following, oft-quoted, prescient remarks of John von Neumann made at the dawn of the computer era---see his collected works edited by Taub [12]---serve to characterize theoretical computer science:

"There exists today a very elaborate system of formal logic, and specifically, of logic applied to mathematics. This is a discipline with many good sides but also serious weaknesses. ……. Everybody who has worked in formal logic will confirm that it is one of the most technically refractory parts of mathematics. The reason for this is that it deals with rigid, all-or-none-concepts, and has very little contact with the continuous concept of the real or the complex number, that is, with mathematical analysis. Yet analysis is the technically most successful and best-elaborated part of mathematics. Thus formal logic, by the nature of its approach, is cut off from the best cultivated portions of mathematics, and forced into the most difficult terrain into combinatorics.

The theory of automata, of the digital, all-or-none type as discussed up to now, is certainly a chapter in formal logic. It would, therefore, seem that it will have to share this unattractive property of formal logic. It will have to be, from the mathematical point of view, combinatorial rather than analytical."

And these observations are echoed by Knuth [8] as follows:

"The most surprising thing to me, in my own experiences with applications of mathematics to computer science, has been the fact that so much of the mathematics has been of a particular discrete type."

These mathematical requisites typically needed to study the algorithms of computer science have been gathered together in Graham, Knuth, and Patashnik [13].

## 2.3 Historical Role of Numerical Analysis within Computer Science

When the discipline of computer science came into existence with the advent of electronic computing during the 1940s, another foundational pillar was knowledge of an array of algorithms, mostly of a numeric nature and associated with the names of famous mathematicians of the past, for example, Euclid, Newton, Euler, Fourier, Gauss, and Cauchy. These classical numerical algorithms were natural candidates for problem-solving on a digital computer. A key application, George Dantzig's linear programming (LP) model and simplex algorithm for efficiently solving LPs, was invented during the mid-1940s and evolved in tandem with advances in electronic computing. In consequence, mathematicians identified with the fields of numerical analysis and optimization played a central role in the creation of computer science departments within universities. Academic departments of mathematics continued to field strong groups in numerical analysis, but many numerical analysts, especially those concerned with solving problems defined over numerical spaces of *finite* dimension, transitioned to newly-created departments of computer science during the 1960s.

A key concern of such numerical analysts in the early days was how to cope with the limitations of floating-point arithmetic. Here the fundamental work of Wilkinson [14], who was a close collaborator of Alan Turing, provided basic guidelines and a number of subtle and beautiful concepts, including backward-error analysis (coupled with perturbation theory), stability of algorithms, ill-conditioning of problems, and so on. But these developments simultaneously served to fix an image of numerical analysis in the eyes of other computer scientists, who sometimes tended to look down on work of the error-analysis variety. Computer science had begun rapidly to move away from computations involving numbers---so-called "number crunching"---and toward the manipulation of (discrete bits of) digital information. It centered increasingly on the "care and feeding of computers," namely, subjects like data structures, programming languages, operating systems, machine organization, theory of computation, artificial intelligence, and so on, and only incidentally on numerical analysis. Over the subsequent decades, numerical analysts began to drift back to departments of mathematics or applied mathematics, sometimes farther afield to operations research or other engineering areas---see the quoted remarks of Nick Trefethen at the beginning of Section 1. In a continuation of these observations, which were published in the *Princeton Companion to Mathematics* (Gowers et al. [1]), he characterizes modern numerical analysis as follows (italics mine):

"In the 1950s and 1960s, the founding fathers of the field [of numerical analysis] discovered that inexact arithmetic can be a source of danger, causing errors in results that "ought" to be right. The source of such problems is *numerical instability*, that is, the amplification of rounding errors from microscopic to macroscopic scale by certain modes of computation. These men, including Von Neumann, Wilkinson, Forsythe, and Henrici, took pains to publicize the risks of careless reliance on computer arithmetic. These risks are very real, but the message was communicated all too successfully, leading to the widespread impression that the main business of numerical analysis is coping with rounding errors. In fact, *the main business of numerical analysis is designing algorithms that converge quickly;* rounding error analysis, while a part of the discussion, is rarely the central issue. If rounding error vanishes, 90% of numerical analysis would remain."

Trefethen further notes that numerical analysis is a discipline that is "built on strong foundations, the mathematical subject of *approximation theory*," and that it has grown into "one of the largest branches of mathematics, the specialty of thousands of researchers who publish in dozens of mathematical journals as well as application journals across the sciences and engineering."

Along similar lines within a classical textbook on numerical analysis, whose title *Analysis of Numerical Methods* serves to encapsulate the subject, Isaacson and Keller [15] make the following observations (italics mine):

"Our opinion is that the analysis of numerical methods is a broad and challenging mathematical activity whose central theme is the effective constructability of various kinds of approximations, …….[and that]  deeper studies of numerical methods would rely heavily on *functional analysis*."

## 2.4 The Foundations of Computational Mathematics (FoCM)

The repatriation of numerical analysts from computer science to mathematics in the 1990s, as just discussed, had an important consequence. The study of fundamental models of computation and complexity developed within theoretical computer science---the so-called "grand unified theory of computation"---leapfrogged back into mathematics, from whence the subject had originated with mathematical logicians of the 1930s. This development was thanks largely to the work of Fields Medalist Stephen Smale and his co-workers. Smale notes the following (quoted from the panel discussion in Renegar, Shub, and Smale [16]):

"A lot of my motivation in spending so much time trying to understand numerical analysis is to help my own ideas about how to define an algorithm. It seems to me that it is important [if one is] to understand the subject of numerical analysis to make a definition of algorithm …….. It is the main object of study of numerical analysis and to have a definition of it so someone can look at all algorithms or a class of algorithms is an important line of understanding."

And he adds:

"…..numerical analysis does not need these things. It doesn't need a model of computation. But on the other hand, I think that [it] will develop. It's going to develop anyway, and it is going to develop probably more in parallel with existing analysis numerical. Numerical analysis will do very fine without it. But in the long run, these ideas from geometry and foundations will give a lot of insights and understanding of numerical analysis."

In a resulting landmark monograph, Blum, Cucker, Shub, and Smale [3] presented a computational model (BCSS) of great generality---abstract machines defined over mathematical rings and fields---and then developed a theory of computational complexity, in particular, over the real and complex numbers. One could summarize this activity, which provides a theoretical foundation for numerical analysis within the field of mathematics, as setting out to bring ""algorithm" under the rubric of "number", the former being appropriately delineated in [3], and the latter being broadly conceived as the branches of mathematical analysis that flowed from the number concept. The foregoing developments led also to the establishment of the Foundations of Computational Mathematics (FoCM) Society [4] to foster such activities.

## 3. Numerical Algorithmic Science & Engineering

### 3.1 Recapitulation

As we have noted in previous sections, numerical analysts within mathematics have been remarkably successful in bringing "algorithm under the rubric of number" and formulating its foundational models. In contrast, computer scientists have been less than successful in achieving the complementary task of bringing "number under the rubric of algorithm." The reason is that numerical computation within computer science has increasingly focused on solving discrete and/or combinatorial numerical problems for which the foundational underpinnings, provided by the classical "grand unified theory" of computation, suffice. But, now recapitulating and briefly enlarging on the introductory discussion of Section 1, it is important to observe that discrete numerical problems arise very often *in conjunction* with continuous, finite-dimensional problems, usually over the real numbers, for example, within mixed integer programs or optimization problems defined over graphs and networks. *The task of building a solid theoretical (and practicable) foundation for numerical computation within computer science as it relates to real numbers, and integrating it with the classical Turing-based models of the subject, was left largely unfinished.* This objective was achieved by numerical analysts during their tenure within CS, but only in a limited, albeit practically important way, through the development and study of the finite-precision, floating-point number model and its associated round-off error analysis.

After numerical analysts migrated back to departments of mathematics, the floating-point model itself has become somewhat of an orphan child. For example, the late Gene Golub makes the following remarks in the panel discussion reported in Renegar, Shub, and Smale [16]:

"I'd like to say something about floating-point arithmetic…. . It is important to know, a few people should know it perhaps. But I don't consider it a part of the mainstream of numerical analysis any longer. Perhaps one needs to know the model. But along with Wilkinson error analysis it isn't in the mainstream of what we call scientific computing today." [Elsewhere in the discussion, he characterizes scientific computing as "a combination of numerical analysis, applied mathematics, and computing."]

### 3.2 Rationale for NAS&E

How does one set about remedying the situation described above? Let us begin by explicitly identifying the discipline within computer science that is the desired *counterpart* of numerical analysis within mathematics:

*Numerical Algorithmic Science and Engineering (NAS&E), or more compactly Numerical Algorithmics, is the theoretical and empirical study and the practical implementation and application of algorithms for solving problems of a numeric nature that are either discrete or continuous over a space of finite dimension and usually defined over the reals, or a combination of the two. This discipline lies within computer science and at its intersection with scientific computing, and it supports the modern modus operandi known as computational science and engineering.*

In contrast, numerical analysis places its emphasis, first and foremost, on continuous problems defied over *function spaces,* i.e., it devotes the lion's share of attention to *infinite-dimensional* numerical problems, for example, partial differential equations, systems of ordinary differential equations, problems of optimal control, etc., and it is undergirded by the mathematical areas of functional analysis and approximation theory (see Section 2.3). Indeed, in our view, the great watershed in numerical computation is between infinite-dimensional numerical problems, on the one hand, and discrete and/or continuous, finite-dimensional problems, on the other. However, it is important to add a qualification that infinite-dimensional problems are generally solved in practice by a reduction to numerical problems of finite dimension, and the latter subject remains an important component of numerical analysis. Thus numerical algorithmics (NAS&E) within computer science as defined above and numerical analysis within mathematics overlap with one another in the domain of continuous, finite-dimensional numerical problem-solving, although often with significant differences in emphasis and the types of algorithmic issues that are addressed, e.g., large-scale, sparse, highly-structured, and so on. In other words, the "watershed" image evoked earlier is more like a flat-topped ridge, with the shared, flat region corresponding to continuous, finite-dimensional problems and the opposing sides corresponding to discrete and infinite-dimensional problems, respectively.

We suggest that specialists in NAS&E be called *numerical algorists*, the counterpart within computer science of numerical analysts within mathematics (see also Nazareth [17]). The word "algorist" has a proud tradition, stretching back to the great Persian mathematician Al-Khowarizm (9th Century, A.D.) from whose name and works both "algorithm" and "algebra" are derived. It is often said that the words "algorithmic thinking" characterize the field of computer science---see, in particular, the 1985 essay "Algorithmic thinking and mathematical thinking" of Donald Knuth that can be found, in expanded form, in Knuth [8]---and, likewise, one could characterize numerical algorithmics (NAS&E) within computer science as "algorithmic thinking applied to discrete and/or continuous finite-dimensional numerical problems."

### 3.3 Foundations of NAS&E

An NAS&E discipline within computer science along lines described above must be based on solid theoretical foundations. To date, the foundational models of computation within CS embrace the following:

1. Classical models of computation developed by mathematical logicians in the 1930s, leading to the so-called "grand unified theory of computation" mentioned earlier (see Moore and Mertens [2]), among which the deterministic and non-deterministic Turing machines models and the random-access machine/random-access stored program (RAM/RASP) models are pre-eminent. For an overview, see Nazareth [18].

2. The modern theory of complexity premised on these models, which has identified problem categories such as P (polynomial-time), NP (non-deterministic, polynomial-time), NP-complete (the most challenging sub-category of NP problems), NP-hard, P-Space, and so on;

3. Cellular automata models developed in Wolfram [19], including universal versions.

4. Randomized models (see, for example, Chapter 10 of Moore and Mertens [2]), and, more generally, quantum models of computation and the universal quantum computer of Deutsch [20] (see, for example, Chapter 15 of Moore and Mertens [2] or Nielsen and Chuang [21] for overviews of this subject).

5. Biological-based, associative models of computation (see Chapter 3 of Churchland and Sejnowski [22]), which have evolved into the modern, multi-layer neural network models for deep learning and the emerging field of data science (see, in particular, the extensive, web-accessible bibliography of Drakopoulos [23]).

Notably missing from the foregoing list are models for *real-number* computation that are well-integrated with the foregoing classical models for discrete, numerical computation. Such models would be the counterpart for NAS&E within computer science of the BCSS model mentioned in Section 2.4, which forms the theoretical foundation for numerical analysis within mathematics.

Let us now briefly survey some progress made to date on this desired CS counterpart:

### 3.3.1 Abstract Real-Number Models

The famed computer scientist Stephen Cook and his co-worker Mark Braverman [24] have promulgated the "bit-model" for scientific computing,

which seeks to remain, as closely as possible, within the Turing tradition. In essence, a real-number function f(x) is computable in their approach if there exists an algorithm which, given a good *rational* approximation to x, finds a good *rational* approximation to f(x). They contrast their "bit-model" with the "algebraic" approach, which the late Joseph Traub [25], another eminent computer scientist, has characterized as follows (italics mine):

"A central dogma of computer science is that the Turing machine is the appropriate abstraction of the digital computer. …. I argue here that physicists [and indeed all scientists] should consider the real-number model of computation as more appropriate and useful for scientific computation. …. The crux of this model is that one can *store and perform arithmetic operations and comparisons on real numbers* exactly and at unit cost."

A "magnitude-based" formalization of Traub's approach, but with *logarithmic costs* for its basic arithmetic operations, is presented in Nazareth [18]. It *re-conceptualizes the floating-point number system* so as to permit computation with real numbers within the standard and well-known RAM/RASP model of theoretical computer science. In the resulting, so-called CD-RAM/RASP model of computation, a real number is defined by a mantissa and an exponent, the former being represented by an *analog* "magnitude," or A-bit, and the latter by a finite, digital sequence of unary (or binary) bits. Arithmetic operations between A-bits are defined *abstractly* by means of geometric-based operations on magnitudes, which are now potentially implementable in practice through the use of HP-memristors---for details, see Nazareth [26].

### 3.3.2 Finite-Precision, Floating-Point Model

The well-known, finite-precision floating-point arithmetic model and its associated error analysis (Wilkinson [14]) can be fully embraced by the foregoing abstract, real-number models and viewed very simply as their *coarsening for practical purposes.* In other words, the floating-point model can be *subsumed* by the "grand unified theory of computer science," once this theory has been suitably broadened to incorporate real-number models of computation.

In its IEEE 754 standardization, a conventional floating-point number is represented by a word length of, say, $n$ bits. This consists of a single sign bit $b$ ; a set of $es$ bits containing an unsigned binary integer, which represents a shifted or biased exponent $e'$, where $0 \leq e' \leq (2^{es} - 1)$, and from which the bias can be removed to obtain the exponent $e$ (a negative, zero, or positive integer); and the remaining $t = (n\text{-}es\text{-}1)$ bits, which represent a normalized mantissa (significand, fraction) $f$ , and where normalized means the leading binary digit of $f$ is 1. Thus a number $x$ is given by $x = (-1)^b f 2^e$. The number of bits within each of the various components that define the finite-precision floating-point number $x$ is fixed. For instance, a 64-bit number has a single sign bit, $es=11$ bits for the exponent, and $t = 52$ bits for the mantissa.

A detailed description of the finite-precision floating-point model and the associated round-off error analysis of its basic arithmetic operations can be found, for example, in Nazareth [27; Chapter 4]. In particular, a key feature of *ideal*, or axiomatic, floating-point arithmetic, upon which the error analysis of algorithms depends, is that

$$\text{fl}\,(x \# y) = (x \# y)\,(1 + \mu),$$

where $x$ and $y$ are any two *representable* floating-point numbers, `#' denotes any one of the four basic arithmetic operations, `fl' denotes the result of that floating-point arithmetic operation, and $|\mu| \leq 2^{-t}$, where the latter quantity is often called a unit in the last place (ulp). Note that $\mu$ will vary depending upon the operation $\#$ and the operands $x$ and $y$, but the *bound* on $|\mu|$ is, in each case, the same number (ulp).

### 3.3.3 Posit-Unum Model

An alternative universal number (unum) arithmetic framework was developed recently in Gustafson [28], [29] and Gustafson and Yanemoto [30], and is summarized in [30] as follows:

"The **unum** (**u**niversal **num**ber) arithmetic framework has several forms. The original "Type I" unum is a superset of IEEE 754 Standard floating-point format...; it uses a "ubit" at the end of the fraction to indicate whether a real number is an exact float or lies in the open interval between floats. While the sign, exponent, and fraction bit fields *take their definition from IEEE 754, the exponent and fraction field lengths vary* automatically, from a single bit up to some maximum set by the user. Type-I unums provide a compact way to express *interval arithmetic*, but their variable length demands extra management. They can duplicate IEEE float behavior, via an explicit rounding function.

The "Type II" unum abandons compatibility with IEEE floats, permitting a clean mathematical design based on the projective reals."

Finally, "Type III" unums, also called *posits,* were proposed by Gustafson and Yanemoto [30] in a radical departure from the floating-point system and its associated IEEE 764 standard. Within a posit representation of $n$ bits, the leading sign bit $b$ is defined as in a floating-point number. A posit has exponent $e$ and fraction $f$ akin to a floating-point number described above, but unlike the IEEE 764 standard, the exponent and fraction parts of a posit do *not* have fixed length. And a posit has an additional category of bits, known as the *regime*, which is also of variable length.

An excellent mathematical account can be found in "*Anatomy of a posit number"* by John D. Cook [31], which is derived, in turn, from Gustafson and Yanemoto [30]. Following the sign bit, the regime has first priority and is defined by a *unary* sequence of length say $m$, comprising either all zeros or all ones, where $m$ can range from a single bit to as many as $n$-1. If the regime is defined by a sequence of 0's then set $k$=-$m$, otherwise, if defined by 1's, set $k$=$m$-1. The remaining bits, if any, up to a maximum allowable number specified by an exogenous parameter *es,* define the exponent, a non-negative integer $e$, s.t. $0 \leq e \leq (2^{es} -1)$. If there are still bits left after the exponent bits, the rest go into the normalized fraction $f$ defined as 1+ the fraction bits interpreted as following a binary point. (e.g, if the fraction bits are 10011, then $f$=1.10011 in binary.)  The posit $x$ is then defined as follows: $x = (-1)^b f 2^{(e+kw)}$ where $w = 2^{es}$. More explanatory detail can be found in Cook [31] who observes that "the primary advantage of posits is the ability to get more precision or dynamic range out of a given number of bits," i.e., posits have "tapered precision [in the sense that] numbers near 1 have more precision, while extremely big numbers and extremely small numbers have less." He notes also that "there's only one zero for posit numbers, unlike IEEE floats that have two kinds of zero, one positive and one negative," and that "there's also only one infinite posit number. For that reason you could say that posits represent projective real numbers rather than extended real numbers. IEEE floats have two kinds of infinities, positive and negative, as well as several kinds of non-numbers."

However, in order to recover some of the beautiful, Wilkinson-type error analysis associated with the finite-precision, floating-point model, in particular, backward error analysis vis-à-vis perturbation theory when seeking to establish the numerical stability of an algorithm vis-à-vis the numerical stability of the problem that it is solving, e.g., an arbitrary triangular system of linear equations with non-zero diagonal elements by forward- or back-substitution (see Nazareth [27], Sec. 4.4), it may be necessary to have a pre-specified, *lower bound* on the number of bits, say $t$, assigned to the normalized mantissa (fraction), as in foregoing subsection 3.3.2, leaving the remaining *(n-t)* bits for the sign, regime, exponent, and possibly additional higher-order bits for the mantissa, as described above. But this negates some of the characteristics of posits mentioned in the previous paragraph and is a subject that requires further exploration. For further discussion of floating-point versus posits, see Greenbaum [32]. More recently, an insightful, in-depth study of posits is given by De Dinechin et al. [33]. They emphasize the serious drawback of posits vis-à-vis floats mentioned above as follows:

"A very useful feature of standard floating-point arithmetic is that, barring underflow/overflow, the relative error due to rounding (and therefore the relative error of all correctly-rounded functions, including the arithmetic operations and the square root) is bounded by a small value $2^{-t}$ where $t$ is the precision of the format. Almost all numerical analysis … is based on this very useful property." [For consistency with the discussion in Section 3.3.2, the quantity `$p$´ in this quotation has been replaced by `$t$´.]

De Dinechin et al. [33] note that "this is no longer true with posits" and, in consequence, "numerical analysis has to be rebuilt from scratch." An ameliorating remedy might be the option of imposing a pre-specified, lower-bound on the number of mantissa (fraction) bits of posits as mentioned above.

### 3.4 NAS&E Content and Organization

The foregoing developments begin to lay a foundation for NAS&E within computer science that conforms to its traditional roots. While such theoretical and practicable foundational models are consequential, the primary focus of NAS&E must always remain on the *scientific* study of the discrete and/or continuous, finite-dimensional algorithms themselves and their *engineered* implementation at all levels. (For a discussion of hierarchical levels of implementation of numerical algorithms, see Nazareth [34].)

Two early and classic works along these lines are Wilkinson [35] and Dantzig [36]. The former is the definitive, path-breaking work on solving systems of linear equations and the algebraic eigenproblem over real (and complex) finite-dimensional spaces. In reference to this work and quoting again from the aforementioned panel discussion in Renegar et al. [16], Beresford Parlett, one of the world's leading experts in matrix computations and numerical analysis, notes the following:

"Even advancing more in time in the field of matrix computations to the sort of bible written by Wilkinson in the 1950s and 60s, he hardly proves a theorem in the book. I've heard people in the audience criticizing the book, because they say it is very inconvenient to use as a reference. The subject really isn't organized very neatly. I think that is actually a true criticism. You sort of have to swallow it whole."

But this gets to the heart of the matter. Wilkinson [35], the book mentioned in the above quotation, is an inspired work that follows a very different paradigm for presenting its algorithmic material, one that has much more in common with the sciences and engineering than it does with mathematics. It represents quintessential NAS&E. And the same can be said of Dantzig's classic, *Linear Programming and Extensions* [36], which was published in the same early period and is another of the crown jewels of NAS&E.

A wide variety of beautiful and powerful algorithms and algorithmic ideas for finite-dimensional equation-solving (linear and nonlinear and often of very large scale) and for discrete and/or continuous, finite-dimensional optimization (linear, non-linear, integer, network, dynamic, stochastic, etc.) have since been discovered and studied in a similar vein; for example, max-flow/min-cut, central path-following, conjugate gradients, quasi-Newton, the duality principle, Nelder-Mead simplex, branch-and-bound, genetic algorithms, bipartite matching/covering, algorithms based on homotopies, Dantzig-Wolfe decomposition of linear programs, and so on, to name but a few. For further detail, see, for example, Moore and Mertens [2] or Nazareth [37], [38].

Turning to the engineering aspect of NAS&E, an outstanding illustration is given by the relevant "recipes" of Press et al. [39]. The engineering of (often highly-complex) implementations of numerical algorithms in a particular computer language and computing environment and the development of appropriate tools that undergird such implementations are an integral part of the NAS&E discipline. For example, Matlab, Python, dialects of C and Fortran, MPI, and so on, are a vital part of the toolkit for implementing algorithms and a numerical algorist should not be merely a competent user of such tools, but should also have some knowledge of what lies "under the hood." Thus, in addition to the mathematical training needed to study discrete and continuous finite-dimensional numerical algorithms, a trained numerical algorist must be cognizant of the techniques that go into the creation of complex data structures, programming languages, compilers and interpreters, operating systems, basic machine design, and so on, subjects at the heart of an education in computer science. The development of implementations and high-quality mathematical software would be a highly respectable activity within NAS&E and academically recognized and rewarded, just as is the case with the creation of non-numeric software within present-day computer science. It is worth noting that the writing of a large piece of mathematical software is as challenging a task as proving a mathematical convergence or rate-of-convergence theorem---harder, perhaps, because the computer is an unrelenting taskmaster.

*Why change horses in midstream by introducing new nomenclature rather than simply retaining the previous term "numerical analysis"?* The answer is that an educational and research curriculum for numerical algorithmics (NAS&E) within computer science differs significantly in character from its numerical analysis counterpart within mathematics. Let us illustrate this point by outlining an introductory course in NAS&E, which could be structured along the following lines:

1. *Algorithmic Foundations*:  1a: Introduction to the formal notion of algorithm and universal machine and the "grand unified theory of computation" (cf. Section 2.2).  1b: Introduction to the fundamental notion of number (cf. Section 2.1). 1c: Bringing "number" under the rubric of "algorithm" via an introduction to "real-number," abstract models of computation and their practical versions, in particular, the standard finite-precision, floating-point and recently-proposed posit-unum arithmetic systems (cf. Section 3.3).

2. *Numerical Algorithmic Science:* 2a: A selection of some of the beautiful algorithms and algorithmic ideas of discrete and finite-dimensional, continuous numerical computing---see, for example, the ones that were listed earlier in this section. 2b: Applications to realistic numerical problems chosen, for instance, from *Numerical Algorithms: Methods for Computer Vision, Machine Learning, and Graphics* by Justin Solomon [40]) 2c: Some exposure to theoretical convergence and rate-of-convergence analysis, but with much greater emphasis being placed on numerical experimentation with algorithms and algorithmic ideas taught in the foregoing item 2a, using, for example, Matlab or Mathematica (see Nazareth [41] for a detailed illustration).

3. *Numerical Algorithmic Engineering:* 3a: A discussion of practical aspects of implementation, for example, data structures, choice of programming language, and so on---see, for example, Nazareth [27], [42], [43].  3b: The CS techniques and "under the hood" design of systems like Matlab or Mathematica used for numerical experimentation in item 2c above.

This is quite different in content from a standard introductory course on numerical analysis within a mathematics or applied mathematics department---see any of the numerous textbooks available, e.g., Kahaner et al. [44]---which generally begins with an introduction to floating-point computer arithmetic, basic roundoff error-analysis and the numerical solution of systems of linear equations via stabilized Gaussian elimination, and then rapidly moves on to other topics such as polynomial interpolation, quadrature, divided-differences, ordinary differential equation-solving, Fourier transforms, and so on, interleaved with optimization-oriented topics like finding the roots of nonlinear equations and the minima of nonlinear functions, in one and several dimensions.

An introductory course in NAS&E within computer science at the undergraduate or graduate level, such as the one outlined above, would typically be followed by a sequence of other courses providing more in-depth coverage. This sequence would take its place alongside course sequences in the standard areas of computer science, for example, data structures, automata theory & formal languages, compiler & interpreter design, operating systems, basic machine design, and so on. The NAS&E course sequence would simply be another available arrow in the computer science quiver!

The mathematical background required for NAS&E within CS would be knowledge of linear algebra and calculus, or basic analysis. One cannot expect a computer science student at the undergraduate or graduate level to be conversant with functional analysis, which is an essential prerequisite for any in-depth study of numerical analysis. However, nothing would prevent a student of NAS&E within computer science from broadening his/her training through course offerings in numerical analysis (and functional analysis as needed) from a mathematics or applied mathematics department. This is analogous to a student of basic machine design within computer science looking to course offerings of an electrical engineering department in seeking to obtain more in-depth instruction in machine hardware.

The creation of small-scale, informally-structured NAS&E research centers within academic departments of computer science can further facilitate the above educational and research objectives. This is discussed in more detail in Nazareth [18, Section 6.6] and is illustrated by two case studies in [18, Chapter 7], one in numerical algorithmic science and the other in numerical algorithmic engineering.

Our hope is that the NAS&E discipline will begin to reoccupy the region within academic departments of computer science that was left vacant following the repatriation of numerical analysts to mathematics (as described in Section 1 within Trefethen's observations quoted from [1]). And thus, looking now to a more distant horizon, it is appropriate to close this section with a second observation, which is taken again from Trefethen's survey of numerical analysis in Gowers et. al. [1]. He has anticipated the emergence of NAS&E as a sub-discipline within the field of computer science as follows (italics mine):

"….the *computer science of numerical analysis is of crucial importance,* and I would like to end with a prediction that emphasizes this side of the subject… . In a world where several algorithms are known for solving every problem, we increasingly find that the most robust computer program will be one that has diverse capabilities at its disposal and deploys them adaptively on the fly. In other words, numerical computation is increasingly deployed in intelligent control loops. I believe this process will continue, just as has happened in many other areas of technology, removing scientists from further details of their computations but offering steadily growing power in exchange. I expect that most of the numerical computer programs of 2050 will be 99% intelligent "wrapper" and just 1% actual "algorithm," if such a distinction makes sense. Hardly anyone will know how they work, and they will be extraordinarily powerful and reliable, and will often deliver results of guaranteed accuracy."

## 4. Concluding Remarks

Numerical algorithmics (NAS&E) as delineated in the previous sections of this article can be viewed as the aforementioned "computer science side of numerical analysis" and, indeed, NAS&E within computer science and numerical analysis within mathematics would serve to *complement* one another and provide *an improved environment for cooperation* (see also Kahan [45]).

NAS&E would be a conduit into computer science of the fundamental concept of "number" drawn from mathematics---see, in particular, Chaitin [46] and other references given there for an extensive treatment of *computable numbers* from both a philosophical and a mathematical perspective.

NAS&E would also serve as a bridge between computer science and the natural sciences and engineering, because new and effective algorithms are often first invented by engineers and scientists to solve particular problems, and only later subjected to a more rigorous mathematical and computational analysis.

And the discipline of NAS&E within the field of computer science would support the modern modus operandi known as computational science & engineering and its so-called "grand challenge" problems of computing, for example, seeking an explanation for the mystery of protein-folding or the intriguing near optimality of the genetic code.

**Acknowledgements**

It is a pleasure to thank Professors Anne Greenbaum and Randy LeVeque and Dr. Eugene Zak for their helpful feedback on this material. They are, of course, not implicated in shortcomings in the views expressed here, these being purely my own.

**Author Information:** John Lawrence (Larry) Nazareth (larrynaz@uw.edu or johnlawrencenazareth@gmail.com) is a professor emeritus at Washington State University (Pullman) and an affiliate professor at the University of Washington (Seattle). He was educated at the University of Cambridge (Trinity College) and the University of California (Berkeley). He lives on Bainbridge Island near Seattle, Washington. Further background information can be found at www.math.wsu.edu/faculty/nazareth .